\documentclass[a4paper]{article}

\usepackage{type1cm}            
\usepackage{newtxtext}          
\usepackage[varvw]{newtxmath}   

\usepackage{makeidx}            
\usepackage{multicol}           
\usepackage[bottom]{footmisc}   
\makeindex                      

\usepackage{amsmath, amsfonts}  
\usepackage{mathtools}          

\usepackage{graphicx}           
\usepackage{tikz}               
\usepackage{subcaption}         

\usepackage{algorithm}          
\usepackage{algpseudocode}      

\usepackage{natbib}             

\usepackage{hyperref}           

\usepackage{xcolor}             
\usepackage[table]{xcolor}      
\usepackage[colorinlistoftodos]{todonotes}  

\usepackage{bbm}                

\usepackage{booktabs}           
\newcommand{\graystyle}[1]{\textcolor{gray}{#1}}  

\newcommand{\firstrevision}[1]{\textcolor{black}{#1}}

\renewcommand{\marginpar}[2][]{}

\usepackage{authblk}
\usepackage{hyperref}

\begin{document}

\title{A Fast Heuristic for Stochastic Steiner Tree Problems}

\author[1,2]{Berend Markhorst\thanks{ORCID: 0000-0002-3383-6862}}
\author[2]{Alessandro Zocca\thanks{ORCID: 0000-0001-6585-4785}}
\author[2]{Joost Berkhout\thanks{ORCID: 0000-0001-5883-9683}}
\author[1,2]{Rob van der Mei\thanks{ORCID: 0000-0002-5685-5310}}

\affil[1]{CWI, Science Park 123, 1098 XG Amsterdam, The Netherlands \par \texttt{berend.markhorst@cwi.nl, r.d.van.der.mei@cwi.nl}}
\affil[2]{VU Amsterdam, De Boelelaan 1105, 1081 HV Amsterdam, The Netherlands \par \texttt{a.zocca@vu.nl, joost.berkhout@vu.nl}}

\date{} 

\maketitle

\begin{abstract}
    Network design under uncertainty arises in countless real‑world settings and can be captured by the Stochastic Steiner Tree Problem (SSTP). Although there are a few approaches specifically tailored to this stochastic optimization problem, there are considerably more state-of-the-art heuristics for its deterministic variant, the Steiner Tree Problem (STP). In this work, we show how to leverage an existing STP heuristic in building a novel method for solving its stochastic variant, the SSTP. This approach is a powerful, yet simple and easy-to-implement way of solving this complex problem. We test our method using benchmark instances from the literature. Numerical results show considerably faster computation times compared to the state-of-the-art, with a gap of approximately 5\%.
\end{abstract}

\vspace{0.5cm}
\section{Introduction}  \label{sec:introduction}
The two-stage Stochastic Steiner Tree Problem (SSTP) and its variants have a wide range of applications in different industries, such as telecommunication network design~\cite{ljubic_stochastic_2017, leitner_decomposition_2018}, wire routing in Very Large-Scale Integration circuits~\cite{hokama_heuristic_2014}, studying protein-interaction networks in bioinformatics~\cite{cheong_solving_2010}, and future-proof ship pipe routing \cite{markhorst_future-proof_2025}. This problem is a generalization of the Steiner Tree Problem (STP)~\cite{ljubic_solving_2021}, in which one must find a minimum-cost subgraph spanning one set of vertices, referred to as \textit{terminals}. Contrary to the STP, in the SSTP, the edge costs and terminals are affected by uncertainty and are revealed after the first stage. We assume there are finitely many possible outcomes for the second-stage parameters and that the probabilities of these scenarios are known.

Approximation algorithms for the SSTP emerged in the 2000s~\cite{nicole_immorilica_costs_2004, gupta_boosted_2004, hutchison_stochastic_2005, swamy_approximation_2006, charikar_stochastic_2007, gupta_lp_2007}, while the 2010s saw a shift toward exact methods, spurred by the 2014 DIMACS Challenge~\cite{dimacs11}. Key SSTP developments in the SSTP literature include a genetic algorithm~\cite{hokama_heuristic_2014}, branch-and-cut techniques~\cite{cheong_solving_2010, ljubic_stochastic_2017}, ILP formulations~\cite{zey_ilp_2016}, and a state-of-the-art decomposition model~\cite{leitner_decomposition_2018}.

Due to the problem's NP-hardness, exact methods like the state-of-the-art algorithm from~\cite{leitner_decomposition_2018} scale poorly with instance size and the number of scenarios, making them unsuitable for real-time or interactive use. In contrast, a heuristic aims to \firstrevision{\marginpar{\textbf{R1:3}}quickly} provide, high-quality solutions, enabling large-scale and adaptive applications in domains such as telecommunications or location planning~\cite{ljubic_solving_2021}. It can also serve as a warm start to speed up exact methods, which can be done with (commercial) solvers, for example Gurobi~\cite{gurobi_2025}.

We now formally introduce the SSTP, similar to~\cite[Definition 1]{leitner_decomposition_2018}. We are given an undirected graph $G=(V,E)$ with the sets of vertices $V$ and edges $E$, first-stage edge costs $c^{(0)}_e \in \mathbb{R}_{+}$ for each edge $e \in E$, and scenario set $S$. \firstrevision{\marginpar{\textbf{R1:4}}Each scenario $s \in S$ occurs with probability $p^{(s)} \in [0,1]$, $\sum_{s \in S} p^{(s)} = 1$, and has (possibly different) second-stage edge costs $c^{(s)}_e \in \mathbb{R}_{+}$ (typically $> c^{(0)}_e$) and its own set of terminals $T^{(s)} \subseteq V$.} We must choose a subset of edges in the first stage, $\bar{E}^{(0)} \subseteq E$ and second-stage, $\bar{E}^{(s)} \subseteq E$, for every scenario $s \in S$, whose union connects all terminals in $T^{(s)}$ and minimizes the \firstrevision{following sum of the first stage costs and the expected second stage costs}:

\begin{equation}
    \sum_{e \in \bar{E}^{(0)}} c^{(0)}_e + \sum_{s \in S} p^{(s)} \sum_{e \in \bar{E}^{(s)}} c^{(s)}_e.
\end{equation}

\begin{figure}[ht]
    \centering
    \begin{minipage}{0.22\linewidth}
        \begin{tikzpicture}[node distance=1cm]
            \node[circle, draw] (B) {$B$};
            \node[circle, draw, right of=B] (C) {C};
            \node[circle, draw, above of=C] (A) {$A$};
            \node[circle, draw, below of=C] (E) {$E$};
            \node[circle, draw, right of=C] (D) {$D$};
        
            \draw[thick] (A) -- (C) node[midway, left, font=\footnotesize] {1.5};
            \draw[thick] (A) -- (D) node[midway, right, font=\footnotesize] {1.5};
            \draw[thick] (B) -- (C) node[pos=0.6, below, font=\footnotesize] {1};
            \draw[thick] (C) -- (D) node[midway, above, font=\footnotesize] {1};
            \draw[thick] (D) -- (E) node[pos=0.7, right, font=\footnotesize] {1};
            \draw[thick] (C) -- (E) node[midway, left, font=\footnotesize] {1.5};
        \end{tikzpicture}
        \subcaption{First stage.}
        \label{fig:example1}
    \end{minipage}
    \hfill
        \begin{minipage}{0.22\linewidth}
        \begin{tikzpicture}[node distance=1cm]
            \node[circle, draw, fill=gray!20] (B) {$B$};
            \node[circle, draw, right of=B] (C) {$C$};
            \node[circle, draw, above of=C] (A) {A};
            \node[circle, draw, below of=C, fill=gray!20] (E) {$E$};
            \node[circle, draw, right of=C, fill=gray!20] (D) {$D$};
        
            \draw[thick] (A) -- (C) node[midway, left, font=\footnotesize] {2};
            \draw[thick] (A) -- (D) node[pos=0.3, right, font=\footnotesize] {3};
            \draw[thick] (B) -- (C) node[midway, below, font=\footnotesize] {1.5};
            \draw[thick] (C) -- (D) node[midway, above, font=\footnotesize] {1.5};
            \draw[thick] (D) -- (E) node[pos=0.7, right, font=\footnotesize] {2};
            \draw[thick] (C) -- (E) node[midway, left, font=\footnotesize] {3};
        \end{tikzpicture}
        \subcaption{2\textsuperscript{nd} stage, $s=1$.}
        \label{fig:example2}
    \end{minipage}
    \hfill
    \begin{minipage}{0.22\linewidth}
        \begin{tikzpicture}[node distance=1cm]
            \node[circle, draw] (B) {$B$};
            \node[circle, draw, right of=B] (C) {$C$};
            \node[circle, draw, above of=C, fill=gray!20] (A) {$A$};
            \node[circle, draw, below of=C, fill=gray!20] (E) {$E$};
            \node[circle, draw, right of=C, fill=gray!20] (D) {$D$};
        
            \draw[thick] (A) -- (C) node[midway, left, font=\footnotesize] {2};
            \draw[thick] (A) -- (D) node[pos=0.3, right, font=\footnotesize] {3.5};
            \draw[thick] (B) -- (C) node[midway, below, font=\footnotesize] {2};
            \draw[thick] (C) -- (D) node[midway, above, font=\footnotesize] {2};
            \draw[thick] (D) -- (E) node[pos=0.7, right, font=\footnotesize] {1.5};
            \draw[thick] (C) -- (E) node[midway, left, font=\footnotesize] {3};
        \end{tikzpicture}
        \subcaption{2\textsuperscript{nd} stage, $s=2$.}
        \label{fig:example3}
    \end{minipage}
    \hfill
    \begin{minipage}{0.22\linewidth}
        \begin{tikzpicture}[node distance=1cm]
            \node[circle, draw] (B) {$B$};
            \node[circle, draw, right of=B] (C) {$C$};
            \node[circle, draw, above of=C] (A) {$A$};
            \node[circle, draw, below of=C] (E) {$E$};
            \node[circle, draw, right of=C] (D) {$D$};
        
            \draw[thick] (A) -- (C);
            \draw[thick, gray!50] (A) -- (D);
            \draw[thick, gray!50, dashed] (B) -- (C);
            \draw[thick, gray!50, dashed] (C) -- (D);
            \draw[thick, gray!50] (D) -- (E);
            \draw[thick] (C) -- (E);
        \end{tikzpicture}
        \subcaption{Optimal solution.}
        \label{fig:example4}
    \end{minipage}
    \caption{Illustrative example of a Stochastic Steiner Tree Problem.}
    \label{fig:small_example}
\end{figure}
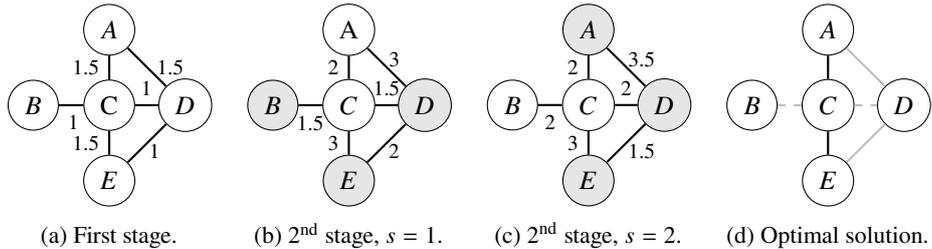

Figure~\ref{fig:small_example} shows a small example of the SSTP with \firstrevision{two equally-likely scenarios}. In Figure~\ref{fig:example1}, we show the graph including the first-stage edge-costs, and in Figures~\ref{fig:example2} and~\ref{fig:example3}, we show the second-stage edge-costs and indicate the terminals with a gray shade for scenario $s=1$ and $s=2$. Figure~\ref{fig:example4} contains the optimal solution, where solid gray edges are chosen in the first stage, and dashed gray edges in the second stage, if necessary.

While few solution methods target the SSTP directly, many exist for its deterministic counterpart, the STP~\cite{ljubic_solving_2021}. This work shows how an STP method can be adapted to solve the SSTP by decomposing it into subproblems, solving them independently, and combining the results. 

The remainder of this work is structured as follows. In Section~\ref{sec:methodology}, we explain our method, after which we provide the setup of our experiments and present the corresponding results in Section~\ref{sec:results}. Finally, we make concluding remarks in Section~\ref{sec:conclusion}.
\section{Methodology} \label{sec:methodology}
We present our heuristic and illustrate the intuition behind it with a small example in Sections~\ref{sec:pseudo-code} and~\ref{sec:small-example}, respectively.

\subsection{Heuristic description} \label{sec:pseudo-code}

\firstrevision{\marginpar{\textbf{R2:5}}Our heuristic in Algorithm~\ref{alg:heuristic}\marginpar{\textbf{R1:5}} works by decomposing the stochastic Steiner tree problem into smaller deterministic subproblems that can be solved quickly using an existing STP heuristic. It first partitions the scenario set into subsets and, for each subset, merges the terminals of all scenarios in that subset into a single deterministic STP. Solving this aggregated problem for first-stage costs produces a tree $\mathcal{E}^{(k)}$ that contains edges which are, loosely speaking, simultaneously cheap in the first stage and useful across many scenarios in that subset. The edges in $\mathcal{E}^{(k)}$ are candidate backbone edges that might be worth buying up front.}

\firstrevision{\marginpar{\textbf{R2:5}}Next, the algorithm considers each scenario separately but restricts itself to the candidate edges identified for its subset. By solving a scenario‑specific STP on this reduced graph, it records which candidate edges are actually used in each scenario. Based on the scenarios in which each edge is used, it compares the one‑time first‑stage cost to the expected second‑stage cost of paying for that edge in each relevant scenario. Edges that are frequently used and relatively expensive in expectation in the second stage are purchased in the first stage; the remaining required edges are then completed per scenario in the second stage. This way, the algorithm exploits the fast deterministic STP heuristic while still capturing the main stochastic trade‑off between upfront investment and recourse. In other words, the heuristic accounts for \marginpar{\textbf{R1:6}}\textit{the value of flexibility}, namely the advantage of a first-stage network that balances expected performance across scenarios and thus can be adapted at relatively low cost. By using multiple scenario subsets, it aims to find edges that are systematically useful across scenarios, guiding more robust first-stage decisions. There is, however, a limit to how much we should aggregate: combining all scenarios creates an artificial problem that encourages excessive upfront spending. Moderate subsets allow the heuristic to find a robust backbone that is neither shortsighted nor overly conservative. \marginpar{\textbf{R1:1}}The appropriate subset size and configuration are problem-dependent and remain a topic for further research.}

The STP heuristic we use has time complexity $O\left(|T^{(s)}||V|^2\right)$ for scenario $s$, see~\cite{kou1981fast}, and therefore, the total time complexity of our heuristic is $O\left(|S| \max_{s \in S} |T^{(s)}||V|^2 \right)$. \firstrevision{\marginpar{\textbf{R1:1}}It is guaranteed that the STP solution of \cite{kou1981fast} is at most a factor $2 \left(1 - \frac{1}{l}\right)$ worse than the optimum, where $l$ is the number of leaves in the optimal tree. Exploring how the STP heuristic choice affects the SSTP solution quality is an interesting topic for future research.}

\subsection{Illustrative example of Algorithm~\ref{alg:heuristic}} \label{sec:small-example}
\firstrevision{\marginpar{\textbf{R2:5}}For illustration, we apply Algorithm~\ref{alg:heuristic} to the example from Section~\ref{sec:introduction} and assume that both scenarios are equiprobable. Take subset size $h=1$ so that we get $K=2$ scenarios subsets. Suppose we get subsets $S_1 = \{1\}$ and $S_2 = \{2\}$. Applying \cite{kou1981fast}, we get $\mathcal{E}^{(1)} = \{(B,C), (C,D), (D,E)\}$ and $\mathcal{E}^{(2)} = \{(A,D), (D,E)\}$ for both subsets, respectively. For this simple example, this leads directly to $\mathcal{F}^{(1)} = \mathcal{E}^{(1)}$ and $\mathcal{F}^{(2)} = \mathcal{E}^{(2)}$.} 

\firstrevision{\marginpar{\textbf{R2:5}}We now complete the algorithm for $e = (A,D)$ as an example. For $e = (A,D)$, it holds that $y^{(1)}_e=0$ and $y^{(2)}_e=1$ and since $c_e^{(0)} = 1.5$ is smaller than $\sum_{s \in S} c^{(s)}_{e} \cdot p^{(s)} \cdot y^{(s)}_{e} = 1.75$, we choose $(A,D)$ in the first stage and thus add it to $\bar{E}^{(0)}$. Informally, this means that buying edge $(A,D)$ in the first stage is cheaper than the expected cost of deferring the purchase, based on the edge's projected usage across the weighted scenarios. Completing the same steps for the remaining edges leads to the optimal solution as shown in Figure~\ref{fig:example4}, i.e., $\bar E^{(0)} = \{(A,D), (D,E)\}$, $\bar E^{(1)} = \{(B,C), (C,D)\}$, and $\bar E^{(2)} = \emptyset$.}

\begin{algorithm}[ht]
\caption{SSTP heuristic.}\label{alg:heuristic}
\begin{algorithmic}[1]
\Statex \firstrevision{\textbf{Input:} Number of scenarios per subset: $h$. Assume $|S|$ is a multiple of $h$ for ease of presentation.}
\Statex \firstrevision{\textbf{Output:} Solution with first-stage edge set $\bar{E}^{(0)}$ and second-stage edge sets $\bar{E}^{(s)}$ for $s \in S$.}
\State Set $y^{(s)}_e \leftarrow 0, \quad \forall s \in S, e \in E$.  \Comment{\graystyle{Equals 1 if edge $e$ is used in scenario $s$, 0 else.}}
\State Set number of subsets $K \leftarrow \frac{|S|}{h}$.
\State Randomly partition scenario set $S$ into $K$ pairwise disjoint subsets $S_1,\ldots,S_K$ with $\bigcup_{k=1}^K S_k = S$.
\ForAll{subsets $k \in \{1,\ldots, K\}$} \Comment{\graystyle{For each partition, find an STP solution.}}
\State Using \cite{kou1981fast}, heuristically solve the STP with terminals $T = \cup_{s \in S_k} T^{(s)}$ based on the first-stage edge costs; collect \firstrevision{the solution} in $\mathcal{E}^{(k)}$.
\EndFor
\ForAll{\firstrevision{scenarios} $s \in S$} \Comment{\graystyle{Asses which edges from $\mathcal{E}^{(k)}$ are used for scenario $s$.}}
\State Set $k$ to the (unique) subset index such that $s \in S_k$.
\State Solve the STP on graph $\bar{G} = (V, \mathcal{E}^{(k)})$ connecting terminals $T^{(s)}$ based on the second-stage edge costs corresponding to scenario $s$; collect \firstrevision{the solution} in $\mathcal{F}^{(s)}$. 
\ForAll{edges $e \in \mathcal{F}^{(s)}$}
\State Update: $y^{(s)}_e \leftarrow 1$.
\EndFor
\EndFor
\ForAll{\firstrevision{edges} $e \in E$} \Comment{\graystyle{Per edge, decide whether it is cheaper to select it in the first stage.}}
\If{$c^{(0)}_{e} \leq \sum_{s \in S} c^{(s)}_{e} \cdot p^{(s)} \cdot y^{(s)}_{e}$}
\State Include $e$ in the first stage solution $\bar{E}^{(0)}$.
\EndIf
\EndFor
\State \textit{The code below is only necessary for obtaining the second-stage solution.}
\ForAll{\firstrevision{scenarios} $s \in S$} \Comment{\graystyle{For each scenario, find a second-stage solution.}}
\State Solve the STP corresponding to scenario $s$ (including terminals and second-stage edge costs), while the edge costs from the edges $\bar{E}^{(0)}$ are set to zero. Store the solution in $\bar{E}^{(s)}$.
\EndFor
\end{algorithmic}
\end{algorithm}
\section{Numerical Experiments} \label{sec:results}
We test numerically the performance of our heuristic. Specifically, we discuss the experimental setup in Section~\ref{sec:experimental_setup} and the corresponding results in Section~\ref{sec:numerical_results}. Finally, we include important remarks in the discussion in Section~\ref{sec:discussion}.

\subsection{Experimental setup} \label{sec:experimental_setup}
Similar to the state-of-the-art method~\cite{leitner_decomposition_2018}, we test the performance of our algorithm on benchmark instances created for the 11th DIMACS Challenge on Steiner Tree problems, which are available on GitHub~\cite{msinnl_sstp}. More specifically, we consider four datasets: K100, P100, LIN, and WRP, whose basic properties are presented in~\cite[Table 1]{leitner_decomposition_2018}. The instances in these datasets range between five and a thousand scenarios.

In preliminary results, we compared five values for $h$ with each other: $\{1,2,5,10,20\}$. This shows that $h = 2$ yields the best results, which is therefore used as the default in the remainder of this work. \firstrevision{\marginpar{\textbf{R1:1}}Due to space limitations, we will not elaborate on the sensitivity analysis on $h$. Yet, this can be a topic for further research. In this work, we consider $h$ primarily as a tunable problem-dependent parameter in our heuristic.}

In this work, we compare our method with the state-of-the-art exact algorithm from~\cite{leitner_decomposition_2018} and the \textit{Buy-None} heuristic from~\cite{hokama_heuristic_2014}. The latter does not acquire any edges in the first stage and waits until the uncertainty is revealed, after which the second-stage solution is computed. We refer to this strategy \firstrevision{as} the ``wait-and-see'' approach and aim to find a balance between the speed of this method and the accuracy of~\cite{leitner_decomposition_2018}.

We \firstrevision{\marginpar{\textbf{R1:8}}implemented} both our algorithm and the wait-and-see approach in Python; our code and data are publicly available on GitHub~\cite{markhorst2025steinerforest}. The experiments are executed single-threaded on a high-performance computing (HPC) cluster with a clock speed of 2.4GHz and 2GB of memory per core. Note that the results from~\cite{leitner_decomposition_2018} are based on a C++ implementation and executed on a single-threaded HPC with a clock speed of 2.5GHz. To present a fair comparison with the state-of-the-art, we decided to report the original results from~\cite{leitner_decomposition_2018} without correcting them with a scaling factor.

\subsection{Numerical results} \label{sec:numerical_results}
We compare our method with~\cite{leitner_decomposition_2018} in terms of solution quality and runtime in Table~\ref{tab:result-table}. The former is quantified by the gap, defined by $\text{Gap} = \frac{UB_H-UB_L}{UB_L}\cdot 100\%$, where $UB_L$ denotes the best upper bound found by~\cite{leitner_decomposition_2018} and provided on GitHub~\cite{msinnl_sstp} and $UB_H$ represents the objective found by our heuristic. In the \firstrevision{\marginpar{\textbf{R1:10}}third} column of Table~\ref{tab:result-table}, we see how often our heuristic obtains a strictly lower gap than wait-and-see. These numbers differ considerably per dataset, which could be explained by the different properties of the underlying graphs. For example, we find a negative correlation between these numbers and the average degree of the graphs. Additionally, we see that the average optimality loss is approximately 5\% for the datasets K100, LIN, and P100, and even around 3\% for the WRP dataset. These solutions are found considerably faster than~\cite{leitner_decomposition_2018} and outperform the wait-and-see heuristic with only a small addition in runtime.

\begin{table}[ht]
    \centering
    \caption{Runtime and gap comparison between our heuristic (H), wait-and-see (W\&S) and~\cite{leitner_decomposition_2018} per dataset.}
    \label{tab:result-table}
    \resizebox{\linewidth}{!}{
    \begin{tabular}{@{}lr>{\columncolor{gray!20}}r>{\columncolor{gray!20}}r>{\columncolor{gray!20}}rrr>{\columncolor{gray!20}}r@{}}
    \toprule
        Dataset & Avg. Degree & Best Gap (\%) & Gap H (\%) & Time H (sec) & Gap W\&S (\%) & Time W\&S (sec) & Time~\cite{leitner_decomposition_2018} (sec) \\ \midrule
        K100 & 7.4 & 14.94 & 4.84 & 3.66 & 5.14 & 1.85 & 17.39 \\
        LIN & 3.3 &70.71 & 4.59 & 97.45 & 5.30 & 56.58 & 1497.35 \\
        P100 & 5.0 & 58.57 & 4.68 & 15.19 & 5.29 & 7.65 & 27.65 \\
        WRP & 3.7 & 61.73 & 2.64 & 90.33 & 2.74 & 49.02 & 2249.34 \\
        \bottomrule
    \end{tabular}}
\end{table}

\subsection{Discussion} \label{sec:discussion}
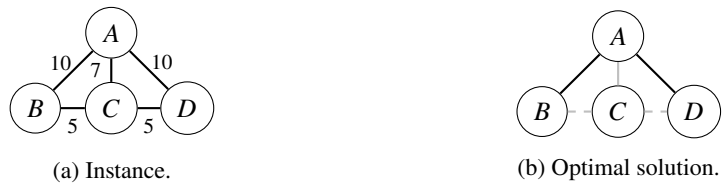
\begin{figure}[ht]
    \centering
    \begin{minipage}{0.45\linewidth}
        \centering
        \begin{tikzpicture}[node distance=1cm]
            \node[circle, draw] (C) {$C$};
            \node[circle, draw, right of=C] (D) {$D$};
            \node[circle, draw, above of=C] (A) {$A$};
            \node[circle, draw, left of=C] (B) {$B$};
        
            \draw[thick] (A) -- (C) node[midway, left, font=\footnotesize] {7};
            \draw[thick] (B) -- (C) node[midway, below, font=\footnotesize] {5};
            \draw[thick] (C) -- (D) node[midway, below, font=\footnotesize] {5};
            \draw[thick] (A) -- (D) node[pos=0.3, right, font=\footnotesize] {10};
            \draw[thick] (A) -- (B) node[pos=0.3, left, font=\footnotesize] {10};
        \end{tikzpicture}
        \subcaption{Instance.}
        \label{fig:drawback_1}
    \end{minipage}
    \hfill
    \begin{minipage}{0.45\linewidth}
        \centering
        \begin{tikzpicture}[node distance=1cm]
            \node[circle, draw] (C) {$C$};
            \node[circle, draw, right of=C] (D) {$D$};
            \node[circle, draw, above of=C] (A) {$A$};
            \node[circle, draw, left of=C] (B) {$B$};
        
            \draw[thick, gray!50] (A) -- (C);
            \draw[thick, gray!50, dashed] (B) -- (C);
            \draw[thick, gray!50, dashed] (C) -- (D);
            \draw[thick] (A) -- (D);
            \draw[thick] (A) -- (B);
        \end{tikzpicture}
        \subcaption{Optimal solution.}
        \label{fig:drawback_2}
    \end{minipage}
    \caption{Drawback of our method: inability to see the value of flexibility. In Figure~\ref{fig:drawback_2}, solid and dotted lines represent the optimal first and second stage solution, respectively.}
    \label{fig:drawback_method}
\end{figure}

\firstrevision{\marginpar{\textbf{R1:11}}At the cost of longer execution times, though still significantly faster than the exact algorithm, the proposed heuristic achieves smaller gaps than the wait-and-see heuristic, especially on the LIN dataset.}
However, these results might heavily depend on the specific features of the instances from~\cite{leitner_decomposition_2018}. For example, the scenarios are almost equiprobable and the ratio between second- and first-stage costs is slightly bigger than one. The latter aids wait-and-see to perform well compared to our method, as selecting everything in the second stage is more affordable than in the first stage. We have enforced this claim with preliminary numerical experiments by increasing the second-stage costs. Based on the results from Table~\ref{tab:result-table}, we observe that our method quickly yields results that outperform wait-and-see. However, a key drawback of our method is that it \firstrevision{still might overlook} the value of flexibility. To illustrate this, we assume that $h=1$, yet the same idea holds for $h>1$. Figure~\ref{fig:drawback_method} shows two equiprobable scenarios: one with terminals $A$ and $B$, the other with $A$ and $D$. Assuming second-stage edge costs are $50\%$ higher, the optimal solution selects edge $(A,C)$ in the first stage and connects $C$ with either $B$ or $D$ in the second. While $(A,C)$ is suboptimal in isolation, it is best when both scenarios are considered. Since our heuristic treats each scenario independently, it misses this flexibility and may underperform in such instances. In this case, it will select nothing in the first stage and either $(A,B)$ or $(A,D)$ in the second stage. The problem in this instance could be addressed by setting $h=2$, yet it is not trivial to prevent in larger problems.

\firstrevision{\marginpar{\textbf{R1:12}}To conclude, the proposed SSTP heuristic is fast and outperforms the wait-and-see heuristic, but is prone to finding sub-optimal solutions due to overlooking the value of flexibility, depending on the choice of $h$ and the specific problem instance being solved.}
\section{Conclusion} \label{sec:conclusion}
We have presented a powerful, yet simple and easy-to-implement method for the SSTP, whose performance is analyzed with runs on benchmark instances from the literature. Using an STP heuristic, we solve the SSTP with fast computation times and a gap of approximately 5\%. Extending this method to the Stochastic Steiner Forest Problem (SSFP) and other stochastic network design problems, e.g.,~\cite{ljubic_stochastic_2017}, has considerable potential to boost research in network design under uncertainty. It is also worth investigating the heuristic's sensitivity to the chosen subset size and the specific STP heuristic used. Another topic for future research is combining our method with scenario reduction, as this could reduce the runtime even more. Finally, adding the Reduced Cost heuristic from~\cite{hokama_heuristic_2014} to our heuristic, possibly in an iterative fashion, could enhance the performance.

\section*{Acknowledgements}
The authors thank Ruurd Buijs and Joris Slootweg for the insightful discussions about this study.

\section*{Competing Interests}
This publication is partly financed by the Dutch Research Council (NWO) with project number \textit{TWM.BL.019.002}.
The authors have no competing interests to declare that are
relevant to the content of this work.

\bibliographystyle{unsrt}
\bibliography{references,extra}

\end{document}